\documentclass[12pt]{amsart}
\usepackage{amscd}
\usepackage{verbatim}
\usepackage{amssymb, amsmath, amsthm, amscd,ifthen}
\usepackage[dvips]{graphics}
\usepackage[cp866]{inputenc}
\usepackage{graphicx}
\usepackage{epsfig}

\usepackage{amsmath,amssymb,amscd,}
\usepackage[all]{xy}

\usepackage{amssymb}

\newtheorem{thm}{Theorem}[section]

\newtheorem{prop}[thm]{Proposition}

\newtheorem{Ex}[thm]{Example}

\newtheorem{lemma}[thm]{Lemma}

\theoremstyle{definition}

\newtheorem{dfn}[thm]{Definition}

\title[Moduli space: a combinatorial
description ]{Moduli space of a planar polygonal linkage: a
combinatorial description}

\author{ Gaiane Panina}

\address{
 St. Petersburg Department of
V.A. Steklov Institute of Mathematics RAS;
St. Petersburg State University, e-mail:gaiane-panina@rambler.ru }

 \keywords{Polygonal linkage, cell complex, CW-complex,  configuration space, moduli space,
permutohedron, cyclic polytope}

\begin{document}
\begin{abstract}
We  describe and study an explicit structure of a regular cell complex
$\mathcal{K}(L)$ on the moduli space $M(L)$ of a planar polygonal
linkage $L$. The combinatorics   is very much related (but not
equal) to the combinatorics of the permutohedron. In particular, the
cells of maximal dimension are labeled by elements of the symmetric
group.
 For example,  if  the moduli
space $M$ is a sphere, the  complex  $\mathcal{K}$ is dual to
the boundary complex of the permutohedron.

The dual complex $\mathcal{K}^*$ is patched of Cartesian products of
permutohedra. It can be
explicitly realized  in the Euclidean space via a
surgery on the permutohedron.

\end{abstract}

 \maketitle \setcounter{section}{0}

\section{Preliminaries and notation}\label{section_preliminaries}

A \textit{polygonal  $n$-linkage} is a sequence of positive numbers
$L=(l_1,\dots ,l_n)$. It should be interpreted as a collection of
rigid bars of lengths $l_i$ joined consecutively in a chain by
revolving joints. We always assume that the triangle inequality
holds, that is, $$\forall j, \ \ \ l_j< \frac{1}{2}\sum_{i=1}^n
l_i$$ which guarantees that the chain of bars can close.

 \textit{A planar configuration} of $L$   is a sequence of points
$$P=(p_1,\dots,p_{n}), \ p_i \in \mathbb{R}^2$$ with
$l_i=|p_i,p_{i+1}|$, and $l_n=|p_n,p_{1}|$. We also call $P$  a
\textit{polygon}.

As follows from the definition, a configuration may have
self-intersections and/or self-overlappings.

\begin{dfn} \textit{The moduli space, or the
configuration space $M(L)$}  is the set  of all configurations of $L$
modulo orientation preserving isometries of $\mathbb{R}^2$.

Equivalently, we can define $M(L)$ as
$$M(L)=\{(u_1,...,u_n) \in (S^1)^n : \sum_{i=1}^n l_iu_i=0\}/SO(2).$$
\end{dfn}

The (second) definition shows that $M(L)$ does not depend
on the ordering of $\{l_1,...,l_n\}$; however, it does depend on the
values of $l_i$.

Throughout the paper (except for the concluding remarks) we assume
that no configuration of $L$ fits a straight line. This assumption
implies that the moduli space $M(L)$ is a closed $(n-3)$-dimensional
manifold  (see \cite{fa}).

The manifold $M(L)$ is already well studied, see \cite{fa,
faS,KapovichMillson}, and many other papers. Explicit descriptions
of $M(L)$ exist for $n=4,5,$ and $6$, see \cite{fa,
KapovichMillson, zvon}. There also exist various results for
polygonal linkages in 3D, see \cite{klya} for example.

The  paper is organized as follows. Section 2 presents an explicit combinatorial description of $M(L)$
as a  regular  cell complex $\mathcal{K}(L)$. In a sense, the starting point of our approach is an elementary version of
Gelfand-Goresky-MacPherson-Serganova idea from \cite{Gelfand}:
they classify the planes (that is, the elements of Grassmanian) by some associated combinatorics.
 The equivalence classes of the planes form \textit{strata} which  may have complicated topology.
   In this paper   we also  classify configurations by their combinatorial types, but here
we are lucky with that all equivalence classes  are topological balls that patch together in a regular cell complex.
The combinatorics of $\mathcal{K}(L)$ is very much related (but not
equal) to the combinatorics of the permutohedron.
In Section 2 we present a number of examples and
give a complete characterization of the possible combinatorics of cells.

In  Section 3 we study
the dual complex $\mathcal{K}^*$ which comes almost automatically with a geometrical realization  in the Euclidean space.
The realization is related to \textit{cyclopermutohedron}  \cite{PanCycloperm},
which is a polytope that encodes cyclically  ordered partitions of a finite set in the same way as the \textit{ permutohedron}
encodes linearly ordered  partitions.

Section 4 sketches the main result of \cite{GalPan}: under a proper setting, a "polygonal linkage" can be replaced by a "simple game" (in the
game-theoretic sense). A simple game cannot be interpreted as a physical object (like bar-and-joint mechanism)
 and therefore has no "configurations".
However,
it is possible to associate with it a cell complex which is proven to be a combinatorial manifold.

Finally, for the sake of completeness,  we discuss in Section 4 the cell complex for the case when the manifold $M(L)$ is singular.

The complex $\mathcal{K}(L)$ already appeared in
\cite{KapovichMillson} in a slight disguise, where it was mentioned
as a "tiling of $M(L)$". Moreover, based on the Deligne-Mostow map,
Kapovich and Millson deduced that
 $\mathcal{K}(L)$ can be  realized as a piecewise linear manifold
in the hyperbolic space.
\bigskip

We start with necessary preliminaries.

\subsection*{Convex configurations}\label{Convexconfigurations}

 A configuration $P$ is \textit{convex }if (1) it is a convex
(piecewise linear) curve, (2) no two consecutive edges are collinear, and (3) the orientation induced by the
numbering goes counterclockwise.

 The set of
all convex configurations we denote by $M_{conv}(L)$. The set
$\overline{M}_{conv}(L)$ is the closure of $M_{conv}(L)$ in $M(L)$.

\begin{lemma}\label{ball} \cite{Kapovich}\begin{enumerate}
    \item The set $M_{conv}(L)$ is an open  subset of $M(L)$  homeomorphic to the open $(n-3)$-dimensional
    ball.
    \item The closure $\overline{M}_{conv}(L)$  is homeomorphic to the closed $(n-3)$-dimensional
    ball.
    \item The interior of $\overline{M}_{conv}(L)$ coincides with $M_{conv}(L)$.
\end{enumerate}
\end{lemma}
Proof. Following   paper \cite{KapovichMillson}, consider  configurations of $n$ (not necessarily all distinct)
points
$p_i$ in the real
projective line $\mathbb{R}P^1$, which we identify with $S^1$.  Each point $p_i$  is assigned the weight
 $l_i$.  The
configuration of (weighted) points is called  {\em
stable} if sum of the weights of coiciding points is
less than half the weight of all points.
%In other words, for a stable configuration, the points $\{p_i\}_{i\in I}$
%may coincide only for a short set $I$.

The group $PSL(2,\mathbb{R})$ naturally acts on the space of
configurations. A remarkable fact is that the quotient  space
of  stable  configurations is exactly
 the space $M(L)$. More detailed, take a stable configuration $\{p_i\}$.
 We interpret the points $p_i$ as unit vectors in $\mathbb{R}^2$.
 In the orbit of the configuration there exists a unique point (up to rotation of $S^1$)
 such that the weighted sum $\sum l_ip_i$ is  zero. Thus each orbit gives a configuration of the linkage $L$.

 A configuration  of points yields a
convex polygon whenever the numbering $(1,...,n)$  goes
counterclockwise. Therefore $M_{conv}(L)$ is identified with the set
of $n$-tuples of  counterclockwise-oriented distinct points $x_i$ in
$S^1=\mathbb{R}P^1$ modulo $PSL(2,\mathbb{R})$. We can omit the
action of the group by assuming that the first three points are $0,\
1 $, and $\infty$. The rest of the points are then given by linear
inequalities
$$
1< x_4< x_5 < ... < x_n < \infty ,
$$  which implies the statement  (1). The statements (2) and (3)  are now straightforward. \qed

\subsection*{Polytopes}\label{pollytopes}

\begin{figure}[h]
\centering
\includegraphics[width=10 cm]{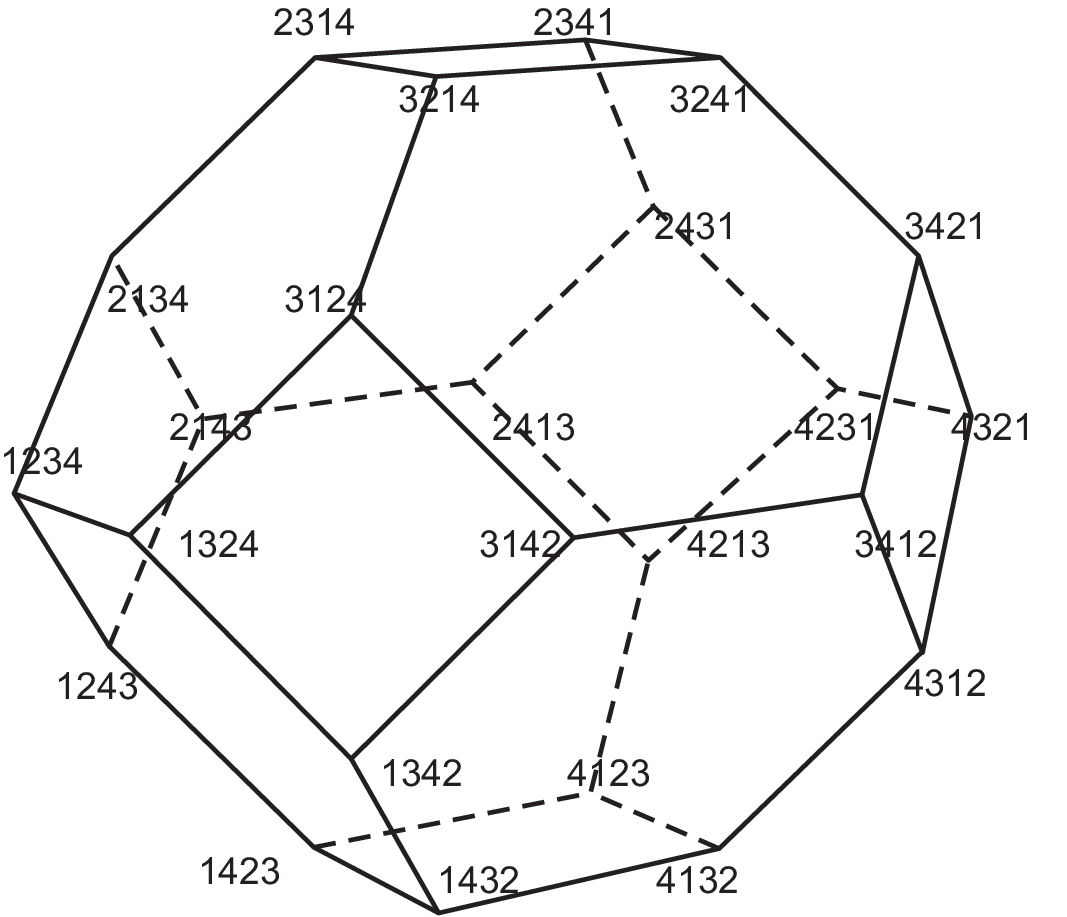}
\caption{Permutohedron $\Pi_4$}\label{permutahedron}
\end{figure}

We shall use the combinatorial structure of the following polytopes:

The \textit{permutohedron}   $\Pi_n$   (see \cite{ziegler} ) is
defined as the convex hull of all points in $\mathbb{R}^n$ that are
obtained by permuting the coordinates of the point $(1,2,...,n)$. It
has the following properties:
\begin{enumerate}
    \item $\Pi_n$ is an $(n-1)$-dimensional polytope.
    \item The $k$-dimensional faces of  $\Pi_n$
are labeled by ordered partitions of the set $\{1,2,...,n\}$  into
$(n-k)$ non-empty parts. In particular, the vertices are labeled by
the elements of the symmetry group $S_n$. The label of a vertex is
obtained by  inverting  the permutation of the coordinates of the
vertex.

    \item A face $F'$ of $\Pi_n$ is contained in a face $F$  iff the label of
$F'$ is finer than the label  of $F$.  Here by a \textit{refinement}
 of an ordered partition $\lambda$ we mean an ordered
refinement  $\lambda'$   whose ordering is inherited from $\lambda$.
For instance, $\{1,3\}\{2,4\}\{5\}$  refines $\{1,3\}\{2,4,5\}$  and
does not refine $\{2,4,5\}\{1,3\}$.
    \item A face of $\Pi_n$ is the Cartesian product of permutohedra of smaller
dimensions.
    \item The permutohedron is a \textit{zonotope}, that is, the Minkowski sum
of line segments.
\item
 The permutohedra $\Pi_1$, $\Pi_2$, and $\Pi_3$ are a one-point polytope,  a
 segment,
 and a regular hexagon respectively. The permutohedron $\Pi_4$ (with its vertices labeled)  is depicted in Fig.
 \ref{permutahedron}.
\end{enumerate}

The \textit{cyclic polytope}  $C(d,n)$   is
 the convex hull of $n$ distinct points $x_1,...,x_n$  on
the moment curve in $\mathbb{R}^d$, see \cite{ziegler}. Its
combinatorics is completely defined by the following property  (\textit{Gale
evenness condition}): a $d$-subset  $F \subset \{x_1,...,x_n\}$ forms
a facet of $C(d,n)$ iff any two elements of $ \{x_1,...,x_n\}
\setminus F$   are separated by an even number of elements from $F$
in the sequence $x_1,...,x_n$.

\section{ The complex $\mathcal{K}(L)$}\label{SectionCW}

\medskip
{\bf Labeling the polygons.} To explain the cell decomposition of
the moduli space, we associate labels to its points.

Assume first that a configuration $P=(p_1,...,p_n)\in M(L)$ has
\textbf{no parallel edges}, that is, no edgevectors
$\overrightarrow{p_ip}_{i+1}$ and $\overrightarrow{p_jp}_{j+1}$ are
parallel and codirected.

Then there exists a unique convex polygon $\overline{P}$  such that

\begin{enumerate}
    \item The edges of $P$ are in one-to-one correspondence with the edges of
    $\overline{P}$. The bijection preserves the directions of the
    vectors.

    \item The  orientations of the edges of
    $\overline{P}$ give the counterclockwise orientation of
    $\overline{P}$.
\end{enumerate}

\begin{figure}[h]
\centering
\includegraphics[width=8 cm]{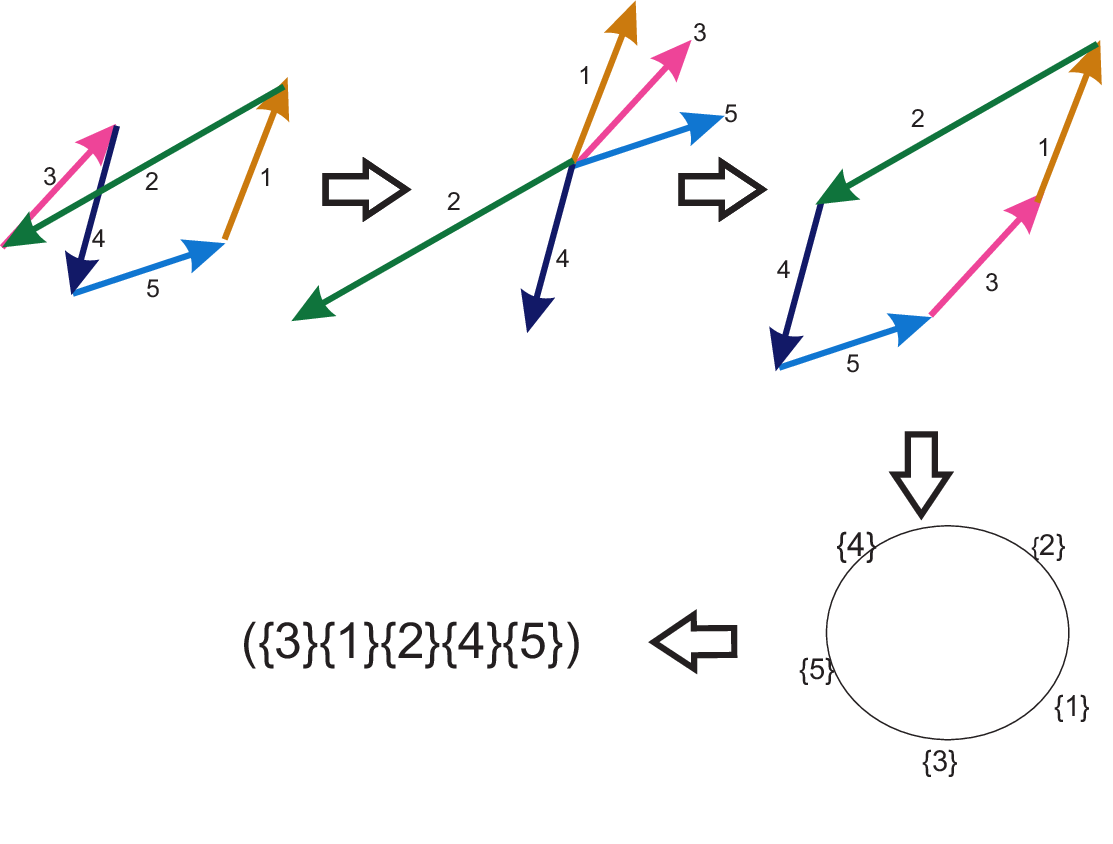}
\caption{Labeling of a polygon with no parallel
edges}\label{rearrangment1}
\end{figure}

In other words, the edges of  $\overline{P}$ are the edges of $P$ coming in the
order  of their slopes (see Fig. \ref{rearrangment1}). Obviously,
$\overline{P} \in M_{conv}(\lambda L)$ for some permutation $\lambda
 \in S_n$. The  permutation is
defined up to the action of the group generated by the cyclic permutation
$(2,3,4,...,n,1)$.   The orbit of a permutation under the action of the group is a \textit{cyclic ordering} on the set $[n]$. Summarizing the above, our construction assigns to $P$ the label
$\lambda(P)$  which is a cyclic ordering  on the set $[n]$.
Equivalently, expecting further discussion on polygons with parallel
edges, we   state that\textit{ a label of a configuration
without parallel edges is a cyclically ordered partition of the set
$[n]=\{1,2,...,n\}$ into $n$ non-empty parts.}

\begin{lemma}\label{lemmaBall1}
Given a cyclically ordered partition $\lambda$ of the set $[n]$ into
$n$ non-empty parts, the subset of $M(L)$ of all polygons labeled by
$\lambda$ is an open $(n-3)$-ball.
\end{lemma}
Proof. The rearranging construction maps the set of  polygons
labeled by $\lambda$ bijectively to  $ M_{conv}(\lambda L)$, which
is a ball by Lemma \ref{ball}.\qed

\begin{dfn}\cite{faS}\label{admiss}
A set $I\subset [n]=\{1,2,...,n\}$ is called \textit{short}, if
$$\sum_{I}^{}l_i <\frac{1}{2} \sum_{i=1}^{n}l_i.$$
\end{dfn}
\begin{dfn}
A partition of the set  $[n]$ is called
\textit{admissible } if all the parts are short.

\end{dfn}

\begin{figure}[h]
\centering
\includegraphics[width=6 cm]{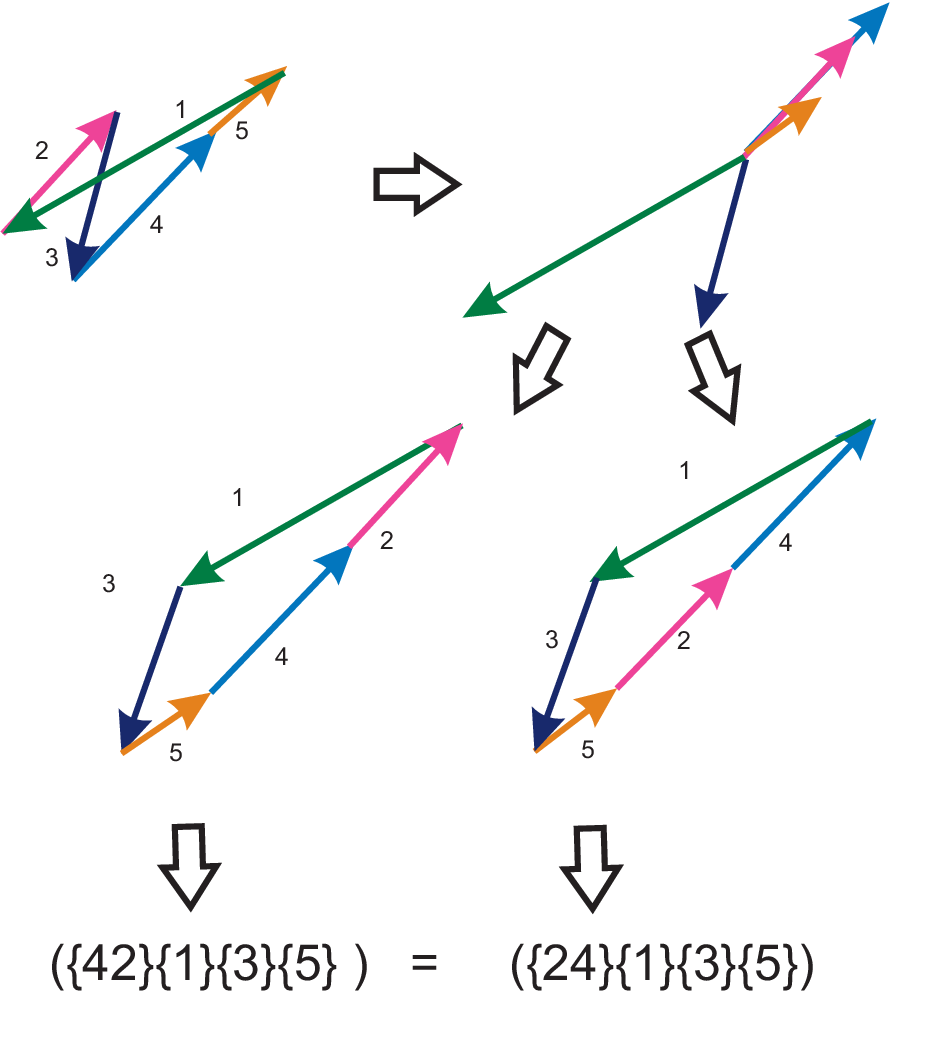}
\caption{Labeling of a polygon with parallel
edges}\label{rearrangment2}
\end{figure}

\bigskip
Assume now that a configuration $P\in M(L)$ \textbf{has parallel
edges}. A permutation which makes $P$ convex is not unique. Indeed,
one can choose any ordering on the set of parallel edges. So in cooking the label,  our
construction puts the indices of parallel edges in one set.

The \textit{label} $\lambda(P)$ assigned to $P$ is\textit{ a
cyclically ordered partition of the set} $[n]$, see Fig.
\ref{rearrangment2} for an example.

\begin{lemma}\label{lemmaBall2}
Given  a cyclically ordered partition $\lambda$ of the set $[n]$
into $k$ non-empty sets, the subset of $M(L)$ of all polygons
labeled by $\lambda$ is either an open $(k-3)$-ball (if $\lambda$ is
an admissible partition), or an empty set (if $\lambda$ is
non-admissible).
\end{lemma}
Proof. We apply Lemma \ref{lemmaBall1} to the $k$-bar linkage with
frozen together edges. Namely, we replace each collection of edges with equal slopes by a single edge. \qed

{\bf A remark on notation.} We write a cyclically ordered partition as a
(linearly ordered) string of sets where the set containing the entry
"$n$" stands on the last position.

We stress once again that the order of the sets matters, whereas there is no ordering inside a set. For
example,
$$(\{1\} \{3 \} \{4,  2, 5,6\})\neq(\{3 \}\{1\}  \{4,  2, 5,6\})= ( \{3 \}\{1\}\{ 2,4, 5,6\}).$$

\begin{dfn}
Two points from $M(L)$ (that is, two configurations) are
\textit{equivalent} if they have one and the same label. Equivalence
classes  of $M(L)$   we  call the \textit{open cells}. The closure
of an open cell in  $M(L)$ is called a \textit{closed cell}. By
above lemmata, all cells are
homeomorphic to balls.

For a cell $C$, either closed or open, its label $\lambda (C)$ is
defined as the label of any interior point of the cell.

\end{dfn}

Before we formulate the main theorem, remind that a CW-complex can
be constructed inductively by defining its skeleta. Once the $(k -
1)$-skeleton is constructed, we attach a collection of closed
$k$-balls $C_i$ by some continuous mappings $\varphi_i$ from their
boundaries $\partial C_i$ to the $(k-1)$-skeleton. For a
\textit{regular } complex, each of  the mappings $\varphi_i$ is
injective, and $\varphi_i$ maps $\partial C_i$ to a subcomplex of
the  $(k-1)$-skeleton, see \cite{Hatcher}.
Regularity of a complex implies that a complex is uniquely defined
by the poset of its cells.
Regularity also guarantees the existence of well-defined barycentric
subdivision and (for PL manifolds) a well-defined dual complex.

\begin{thm}\label{MainThm}
The above described collection of  cells yields  a structure of a
regular CW-complex $\mathcal{K}(L)$ on the moduli space $M(L)$. Its
complete combinatorial description reads as follows:
\begin{enumerate}
    \item  $k$-cells of the complex $\mathcal{K}(L)$ are labeled by cyclically ordered admissible partitions of
 the set  $[n]$  into $(k+3)$ non-empty
parts.

    \item A closed cell $C$ belongs to the boundary of some other closed cell
    $C'$  iff  the partition  $\lambda(C')$ is finer than
    $\lambda(C)$.

\end{enumerate}
\end{thm}
Proof:  The open cells are balls by Lemmata \ref{lemmaBall1} and \ref{lemmaBall2}. The regularity of the complex follows from
Lemma \ref{ball}, (3).\qed

\bigskip

 For the complex $\mathcal{K}(L)$ we immediately have:

\begin{prop}

\begin{enumerate}\item The facets of the complex (that is, the cells of
maximal dimension  $n-3$) are labeled by cyclic orderings on the set
$[n]$.
\item The vertices of the complex are labeled by cyclically ordered admissible partitions of
 the set  $[n]$  into three non-empty
parts.

In other words, they correspond to all possible (oriented) triangles composed
of segments of lengths $l_1,...,l_n$.

\item The  vertex figure   of any vertex $v$ of the complex $\mathcal{K}(L)$ is  combinatorially dual to the
 Cartesian product of three permutohedra.

 More precisely, the label $\lambda(v)$ consists of three parts. If the three parts  have $k$, $l$, and $m$ elements respectively,
 then the vertex figure of $v$ is  combinatorially dual to $\Pi_k\times \Pi_l \times \Pi_m$.
 \item The  face figure of any $k$-dimensional face is  combinatorially dual to the
 Cartesian product of $(k+3)$ permutohedra. (Some of these
 permutohedra can be $\Pi_1$, and thus degenerate to a point.)
\end{enumerate}
\end{prop}

The proof follows directly from the above construction. \qed

\begin{figure}[h]
\centering
\includegraphics[width=8 cm]{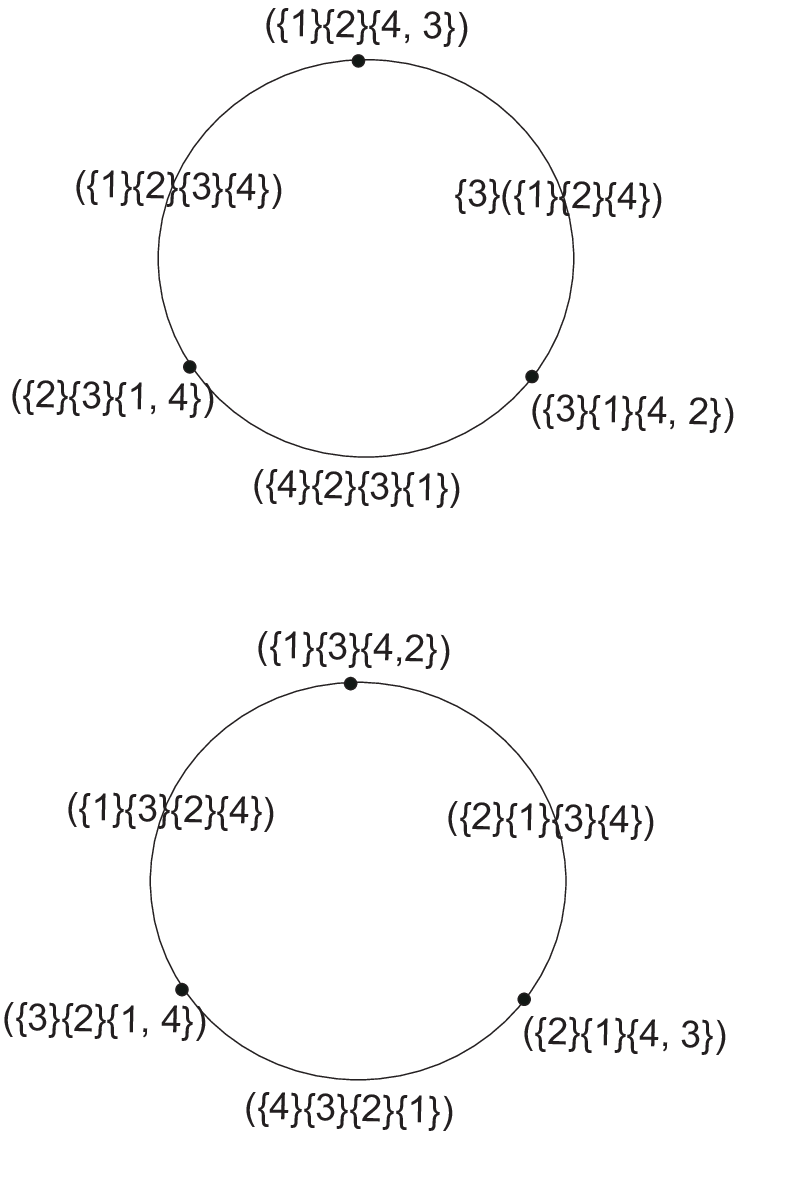}
\caption{$\mathcal{K}(L)$ for the 4-gonal linkage
$(1,1,1,1/2)$}\label{4gon1}
\end{figure}

\begin{Ex}  Let $n=4; \ \ l_1=l_2=l_3=1,\ l_4=1/2.$
The moduli space $M(L)$ is known to be a disjoint union of two
circles, see \cite{fa}. The cell complex $\mathcal{K}(L)$ is
depicted in Fig. \ref{4gon1}.
\end{Ex}

\begin{Ex}\label{ExPErmuto}
Assume that
$$\forall i \ \ l_n+l_i>\sum_{ n\neq j\neq i} l_j.$$
In this case the moduli space $M(L)$ is an $(n-3)$-sphere, see
\cite{fa}, and
the complex $\mathcal{K}(L)$ is dual to the boundary complex of the
permutohedron $\Pi_{n-1}$.
\end{Ex}
Proof. Indeed, each admissible partition is of the type
$$(*, \{n\}),$$ where "$*$" is any linearly ordered partition of $[n-1]$ in at least two parts. This means
that the facets of $\mathcal{K}(L)$ are in a natural bijection with
the vertices of $\Pi_{n-1}$. It remains to observe that the patching
rules for $\mathcal{K}(L)$ are exactly dual to those  of
the permutohedron.\qed

In a regular complex, the boundary of each cell is a combinatorial sphere, so it makes sense to speak of \textit{combinatorics of a cell}.
Let us look what types of combinatorics do we encounter in complexes $\mathcal{K}(L)$ for different linkages $L$.

\begin{Ex}Let $n=5$, $L=(1,1,1,1,1)$. Then $\mathcal{K}(L)$ is a
surface of genus four patched of 24 pentagons. Each vertex has $4$
incident edges. The complex is \textit{flag-transitive},
which means that any combinatorial equivalence of any two pentagons
extends to an automorphism of the entire complex.
\end{Ex}

\begin{Ex}\label{Exequilateral}Let $n=2k+1$, $L=(1,1,...,1)$. Then $\mathcal{K}(L)$ is
 patched of $(2k)!$ copies of duals to the cyclic polytope $C(n-3,n)$.

 However, unlike the
 previous example, the complex is
not completely transitive, just \textit{facet-transitive}: for every
two facets there exists an automorphism of $\mathcal{K}(L)$ mapping one
facet to the other.
\end{Ex}
Proof. Fix a facet $C$ of $\mathcal{K}(L)$. Without loss of generity
we may assume that its label is $(\{1\}\{2\}\{3\}\{4\}...\{n\})$.
Consider the following "starlike" bijection $\varphi$ which maps the
vertices $x_1,...,x_n$ of the cyclic polytope $C(n-3,n)$ to facets
of the cell $C$:
$$\varphi(x_{2i+1})=(\{1\}\{2\}\{i+1,i+2\}\{3\}\{4\}...\{n\}),$$
$$\varphi(x_{2i})=(\{1\}\{2\}\{k+i-1,k+i\}\{3\}\{4\}...\{n\}).$$
Informally, the  defining  rule of $\varphi$ is the way of drawing  a polygonal star (say, a pentagram).
 It is easy to check that $\varphi$ yields a combinatorial duality.\qed

\begin{prop}\label{Propface}
\begin{enumerate}
    \item {\bf Faces of $\mathcal{K}(L)$ are combinatorially equivalent to convex polytopes.} Let  $C$ be a closed cell of $\mathcal{K}(L)$  for some polygonal linkage $L$.
 The boundary complex of $C$  is combinatorially equivalent to a simple $k$-polytope with at most $k+3$ facets.
     Moreover, there exists  some even $D\in \mathbb{N}$
    such that the boundary complex of $C$  is combinatorially
equivalent to a face of the dual to the cyclic polytope $C(D,D+3)$.

    \item {\bf Universality property.} Conversely, any simple $k$-dimensional polytope  $K$ with
    at most $k+3$ facets arises in this way.
    That is, there exist a number $n$, an $n$-linkage
    $L$, and a cell $C$ of the complex $\mathcal{K}(L)$ such that
    the boundary complex of $C$ is combinatorially equivalent to
    $K$.
\end{enumerate}
\end{prop}
Proof. (1) We may assume that all $l_i$ are integers, and that their
sum $D+3=\sum l_i$  is odd. Indeed, neither a small perturbation nor a scaling changes the combinatorics of the complex.
 The space $M(L)$ embeds in a natural
way in the moduli space of the equilateral polygon with $D+3$ edges
$M(\underbrace{1,1,...,1}_{D+3}).$  The embedding maps a polygon
with edgelengths $l_1,...,l_n$ to the equilateral polygon which
represents the same curve, that is, with first $l_1$ edges parallel,
next $l_2$ edges parallel, etc. The embedding respects the structure
of  cell complexes, and therefore, realizes  $\mathcal{K}(L)$ as a
subcomplex of the complex
$\mathcal{K}(\underbrace{1,1,...,1}_{D+3})$, whose facets are
combinatorial cyclic polytopes (see Example \ref{Exequilateral}).

(2)  Assume that a simple $k$-dimensional polytope  $K$ has $k+3$ facets. Then the dual polytope $K^*$
has $k+3$ vertices. We shall prove that every simplicial
$k$-polytope with at most $k+3$ vertices is a face figure of the
cyclic polytope $C(D,D+3)$ for some even $D$. The Gale diagram of
$K^*$ (see \cite{ziegler}) is a one-dimensional configuration of
distinct black and white points. Remind that the Gale diagram of
$C(D,D+3)$ is the alternating configuration of distinct black and
white points in the straight line. Being translated to the Gale
diagram's language, the  statement we need reads as "any
configuration of distinct black and white points in the straight
line can be completed to an alternating configuration of distinct
black and white points", which is obvious.

If $K$ has less than $k+3$ facets, the proof is even simpler.\qed

\section{The dual complex $\mathcal{K}^*(L)$. Surgery on the permutohedron}\label{generic}

\begin{thm}\label{ThmPL} The dual cell complex $\mathcal{K}^*(L)$ carries
a natural structure of a polyhedron.
\end{thm}
Proof. The cells of the dual complex $\mathcal{K}^*$ are the
duals to the face figures of $\mathcal{K}(L)$. By Theorem
\ref{MainThm}, the latter are combinatorially equivalent to
Cartesian products of permutohedra. To realize
$\mathcal{K}^*$ as a polyhedron, for each facet of
$\mathcal{K}^*$ we take the Cartesian product of three standard
permutohedra. Their faces that are identified via isometries.\qed

\bigskip

We describe below a realization of  $\mathcal{K}^*$  in the Euclidean space  $\mathbb{R}^{n-2}$.
For this, we need a preliminary construction which is the subject of paper \cite{PanCycloperm}.
The construction involves the theory of \textit{ virtual polytopes} developed originally in \cite{PuKho}, and some related technique.
For the very first orientation we recommend the reader  just to trust that
there exists a well-defined \textit{Minkowski subtraction} of convex polytopes, and that Minkowski differences have a well-defined facial structure.
For more  details, we refer  to the above mentioned paper.

\subsection*{Cyclopermutohedron}
For a fixed   number $n\geq3$, we define the following  regular cell complex ${CP}_{n}$ by listing all the closed
cells  together with the incidence relations.

 \begin{itemize}
   \item For  $k=0,...,n-3$, the   $k$-dimensional cells ($k$-cells, for short) of the complex ${CP}_{n}$ are labeled by (all possible)  cyclically ordered
partitions of the set $[n]$  into $(n-k)$  non-empty
parts.
   \item A (closed) cell $F$ contains a cell $F'$ whenever the  label of $F'$ refines
 the  label of $F$.
 \end{itemize}

 The complex ${CP}_{n}$ cannot be represented by a convex polytope, since   it is not a combinatorial sphere (not even a combinatorial manifold).
 However, it can be represented by some  \textit{virtual polytope}
 which we call {cyclopermutohedron}  $\mathcal{CP}_{n}$.

Here is the construction of cyclopermutohedron:

Assuming that $\{e_i\}$  are standard basic vectors in $\mathbb{R}^{n-1}$, define the points

$$\begin{array}{ccccccccc}
    R_i=\sum_{j=1}^{n-1} (e_j-e_i)=(-1, & ... & -1, & n-2, & -1, & ... & -1, & -1, &-1, )\in \mathbb{R}^{n-1},\\
     &  & & \ i &  &  &  &  &
  \end{array}
$$
and  the following two families of line segments:
$$q_{ij}=\left[e_i,e_j\right], \ \ \  i<j$$ and $$  r_i=\left[0,R_{i} \right].$$

We also need the point $S=\left(1,1,...,1\right)\in \mathbb{R}^{n-1}$.

\begin{dfn} The  \textit{cyclopermutohedron} is a virtual polytope defined as    the weighted Minkowski sum
of line segments:
$$ \mathcal{CP}_{n}:=   \sum_{i< j} q_{ij} + S- \sum_{i=1}^{n-1}  r_i.$$
\end{dfn}
\begin{thm}\label{thmCycloperm}\cite{PanCycloperm}
 The poset of  (proper) faces of $ \mathcal{CP}_{n}$ is combinatorially  isomorphic to the complex $CP_{n}$.
\end{thm}

\textbf{Remark.}  The  sum
$ S+  \sum_{i< j} q_{ij}$  equals the standard permutohedron.

\medskip

In an oversimplified way, the cyclopermutohedron
$\mathcal{CP}_{n}$  can be visualized as the permutohedron $\Pi_{n-1}$
"with diagonals". This means that all the  proper faces of  $\Pi_{n-1}$
are also faces of  $\mathcal{CP}_{n}$. However,
$\mathcal{CP}_{n}$ has some extra faces in comparison with
$\Pi_{n-1}$.

\medskip

For any $n$-linkage $L$, the complex
 $\mathcal{K}^*(L)$ automatically embeds in ${CP}_{n}$, and therefore embeds in the face complex of $\mathcal{CP}_{n}$.
The embedding goes as follows.
Take the permutohedron $\Pi_{n-1}\subset
\mathbb{R}^{n-1}$,
 assuming (as usual) that the faces  of $\Pi_{n-1}$ are labeled by
ordered partitions on the set $[n-1]$. In particular, the
vertices of $\Pi_{n-1}$ are labeled by permutations of the set
$[n-1]$. We introduce the following bijection between the
vertex sets
$$\psi: Vert(\mathcal{K}^*)\rightarrow Vert(\Pi_{n-1}).$$
Given a vertex of $\mathcal{K}^*$ whose label $\lambda$ is a
cyclically ordered set $[n]$, the mapping $\psi$ sends it  to the vertex
of $\Pi_{n-1}$ by cutting  $\lambda$ at the position of "$\{n\}$" and
omitting "$\{n\}$" from the label.

Thus, the vertices of $\mathcal{K}^*$ are geometrically realized by vertices of the permutohedron.
Next, we realize the cells of the complex:
take a cell $C$  and patch  the  face
of the cyclopermutohedron which corresponds to $C$ by Theorem \ref{thmCycloperm}.

\bigskip

This construction can be reformulated as the following \textbf{surgery algorithm: }

\begin{enumerate}
    \item Start with the complex $\mathcal{K}^*(L)$ and the boundary complex of the permutohedron $\Pi_{n-1}$.
 Realize the vertices of $\mathcal{K}^*$ as the vertices of
$\Pi_{n-1}$ via the above described mapping $\psi$.
    \item For every  face $F$ of $\Pi_{n-1}$ do the following. The face is labeled by some $\lambda$, which is a linearly ordered partition
    of $\{1,...,n-1\}$.
    If the partition is  admissible (that is,  all the parts are short), keep the face $F$ and assign to it
    the label $(\lambda, \{n\})$.
 If the partition is  not admissible, remove the face $F$ from the  complex.

    This step gives a realization of all the cells of $\mathcal{K}^*$
 whose label  contains the one-element  set $\{n\}$.

    \item Take all the cells $C$ of $\mathcal{K}^*$ such that the part of $\lambda(C)$
    containing $n$  has more than one element.  Patch in the corresponding face of  the cyclopermutohedron, which up to a translation equals

    $$\sum  q_{ij} - \sum r_i, $$
    where
        the first (Minkowski) sum extends over all $i<j<n$ such that
    $i$ and $j$ belong to one and the same part of the partition $\lambda(C)$, and
         the second sum extends over all $i<n$ such that
    $i$ and $n$ belong to one and the same part of the partition $\lambda(C)$.
This  is a virtual polytope with the vertex set
$\psi (Vert(C))$.

\end{enumerate}

\begin{Ex} Let $L$ be as in Example \ref{ExPErmuto}.
 The above
described surgery leaves the permutohedron as it is. That is, all
the faces of $\ \Pi_{n-1}$ survive on the second step of the surgery
algorithm, and nothing is added on the third step.

Important is that the "long" edge is the last one. Otherwise we
would get another surgery, but, of course, an isomorphic
combinatorics.
\end{Ex}

\begin{Ex} Let $n=5$; $l_1=1,2;\ l_2=1;\ l_3=1;\ l_4=0,8;\ l_5=2,2$.
 The surgery algorithm  starts with the permutohedron $\Pi_4$  (see
Fig. \ref{torus}). The two shadowed faces are labeled by
$(\{123\}\{4\})$  and $(\{4\}\{123\})$. Since the partitions
$(\{123\}\{4\}\{5\})$  and $(\{4\}\{123\}\{5\})$ are non-admissible,
according to the algorithm, the faces are removed. All other faces
of the permutohedron survive the surgery. Step 3 gives six new
"diagonal" rectangular faces. They correspond to the cells labeled
by $(\{1\}\{2\}\{3\}\{45\})$, $(\{1\}\{3\}\{2\}\{45\})$,
$(\{2\}\{1\}\{3\}\{45\})$, $(\{2\}\{3\}\{1\}\{45\})$,
$(\{3\}\{1\}\{2\}\{45\})$, and $(\{3\}\{2\}\{1\}\{45\})$.
\end{Ex}

\begin{figure}[h]\label{torus}
\centering
\includegraphics[width=12 cm]{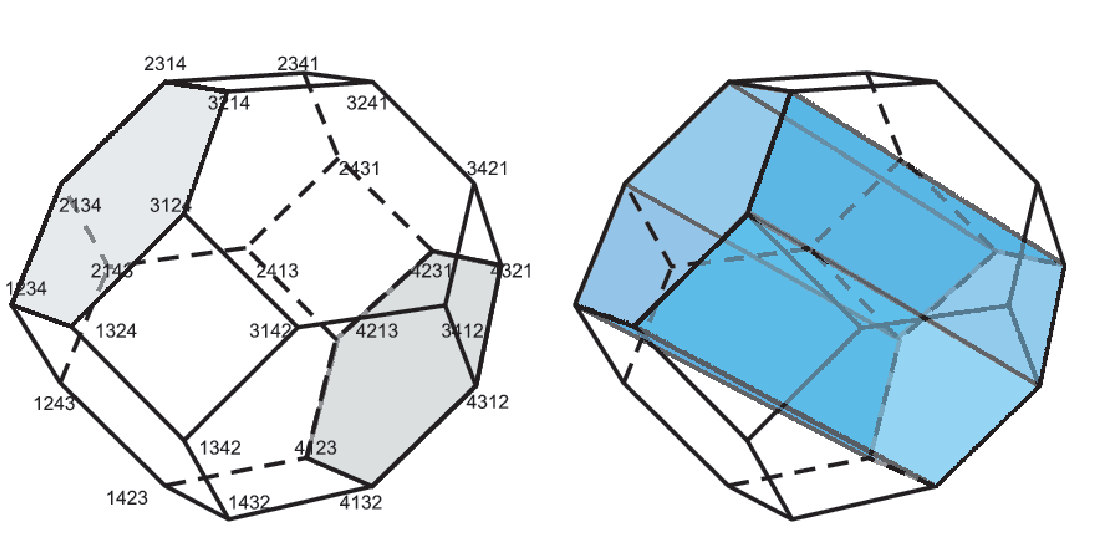}
\caption{The complex $\mathcal{K}^*(L)$ for the $5$-linkage \newline
$L=(1,2;\ 1;\ 1;\ 0,8;\ 2,2)$. We remove from the permutohedron the
grey facets and patch in the blue cylinder.}\label{torus}
\end{figure}

\begin{figure}[h]
\centering
\includegraphics[width=10 cm]{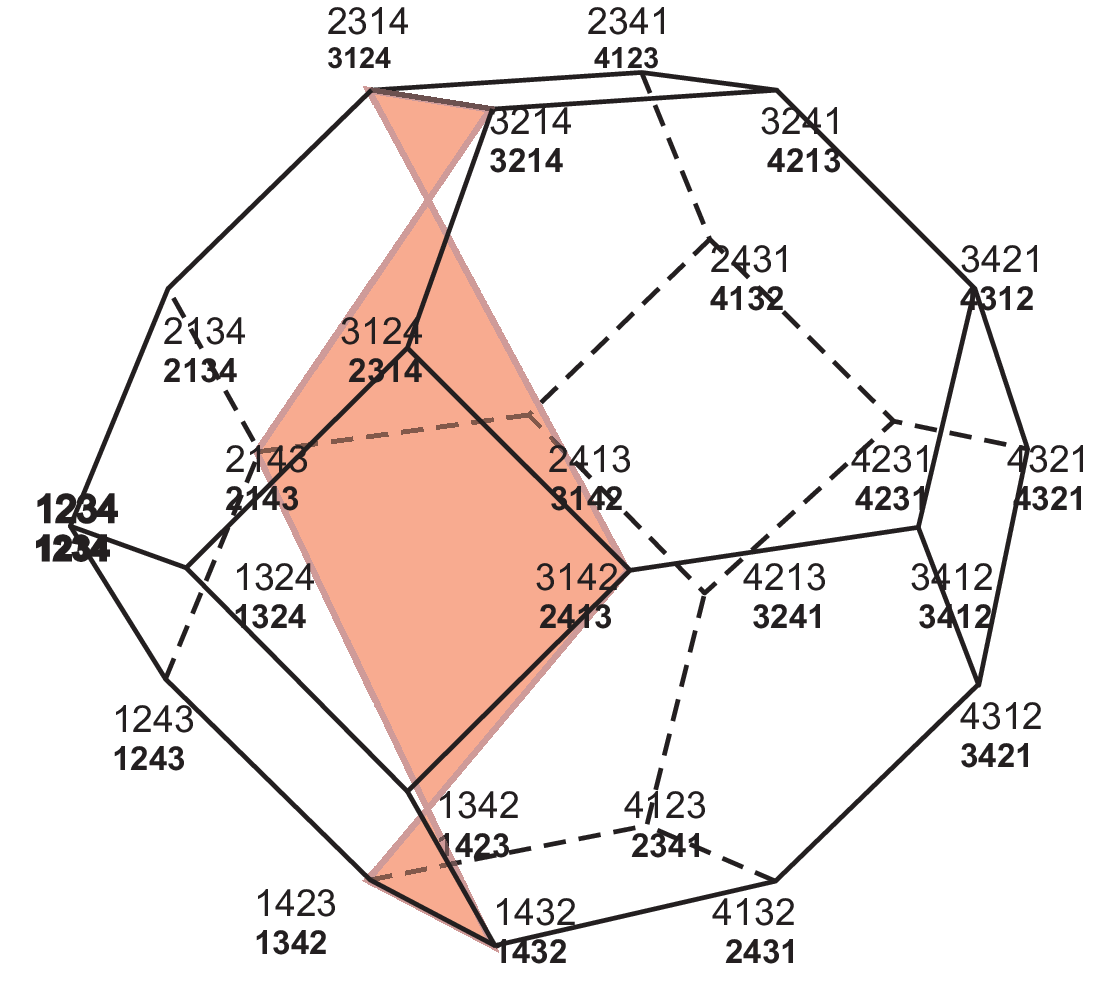}
\caption{ A "diagonal" face }\label{diag}
\end{figure}

\begin{Ex}\label{ExHexaFace} Let $n=5$, $L=(3,\ 1,\ 1,\ 4,\  4)$.
Figure \ref{diag} presents the permutohedron, the labels of the
vertices, and the coordinates of the vertices (in bold). We also
depict the hexagonal face labeled by $(\{1\}\{4\}\{235\})$. It is
the Minkowski sum of two negatively weighted and one positively weighted segments.
\end{Ex}

For more examples of the surgery see \cite{Gor}, where I. Gorodetskaya presented  the surgery for all types of  $5$-linkages.

\newpage

\section{Concluding remarks}\label{concluding}
%\subsection*{Singular configuration space}\label{singular}

The  construction of $\mathcal{K}$ and $\mathcal{K}^*$  suggests some further natural discussions sketched briefly
in this section.

\subsection*{Quasilinkages, simple games,  Alexander self-dual complexes, and associated manifolds}

An elementary observation is that the complex $\mathcal{K}(L)$
depends only on the collection of admissible partitions. In turn, these are defined by the collection of  short sets.
This suggests the following generalization, which is described in details in \cite{GalPan}, and which
we sketch
 very briefly now.

\begin{dfn}
A family $\mathcal{F}$ of subsets of $[n]$ is called {\it a quasilinkage}, if it satisfies the following properties:
\begin{enumerate}
\item $\mathcal{F}$ contains all singletons:  for any $i \in [n]$, $\{i\} \in \mathcal{F}$.
\item Monotonicity: if $S\in \mathcal{F}$, and $T\subset S$ then $T\in \mathcal{F}$.
\item Strong complementarity: if $S\in \mathcal{F}$ then $([n]\setminus S)\notin \mathcal{F}$ , and, conversely, if $S\notin \mathcal{F}$, then $([n]\setminus S)\in \mathcal{F}$.
\end{enumerate}
\end{dfn}

The proposed notion exists in the literature; yet in  completely
different frameworks. It appeared  as ``simple game with constant
sum'' in game theory,  as ``strongly
complementary simplicial complex'', and as " Alexander self-dual simplicial complex".

Being motivated by polygonal linkages, we
call any $S\in \mathcal{F}$  \textit{ a
short set}, and any $S\notin \mathcal{F}$ {\it a long set}.

 Each polygonal linkage $L$ yields a collection of short sets, and therefore, is a quasilinkage. The converse is not true:
 there  exist many quasilinkages that cannot be represented by length assignments.

We  associate with a quasilinkage
$\mathcal{F}$  a cell complex  $\mathcal{K}(\mathcal{F})$ by applying the rules from
Theorem \ref{MainThm}. In \cite{GalPan} it is proven that the complex is a (combinatorial) manifold of
dimension $(n-3)$ which is   locally isomorphic to
$\mathcal{K}(L)$ for some  linkage $L$ (however, $L$ depends on the
location, and  there may be no  linkage associated to the entire
complex).

\subsection*{Cell decomposition for singular configuration
spaces} A similar cell complex exists also for  singular
configuration spaces, that is, for the case when $L$ has
configurations that fit in a straight line.

\begin{dfn}For a singular case,
a partition of  $L=(l_1,\dots ,l_n)$ is called \textit{admissible }
if one of the two conditions holds:
\begin{enumerate}
    \item The number of the parts is greater than $2$, and the total length of any part is strictly
    greater than the total length of the rest.
    \item The number of parts equals $2$, and the lengths of the
    parts are equal.

\end{enumerate}
\end{dfn}

The combinatorics of the complex $\mathcal{K}(L)$ is literally the
same as in Theorem \ref{MainThm} except for the following items:

\begin{enumerate}
\item Non-singular vertices are  labeled by admissible partitions with
exactly three parts.
\item Singular vertices are  labeled by admissible partitions with
exactly two parts.
 \item  Assume that a singular vertex $v$ of $\mathcal{K}(L)$ corresponds to an ordered
 partition of $\{1,2,...,n\}$ into two non-empty parts, say, with
 $k$ and $l$ elements.
 Then the vertex figure of $v$ is  combinatorially equivalent
 to the cone over  $(\partial \Pi_k\times \partial \Pi_l)^* $.
\end{enumerate}

\subsection*{Acknowledgements.} The research was supported by RFBR grant 15-01-02021.
 I'm grateful to Nikolai Mnev for
inspiring conversations. I'm also indebted to Misha Kapovich for
delivering me the proof of Lemma \ref{ball}.


\begin{thebibliography}{99}




\bibitem{Ewa96}
G. Ewald,  Combinatorial convexity and algebraic geometry.
 Springer Verlag, 1996.

\bibitem{GalPan} P. Galashin, G. Panina, Manifolds associated to simple games, J. of Knot Theory and Its Ramifications,   25, No. 12 (2016).

\bibitem{fa} M. Farber, Invitation to Topological Robotics. Zuerich Lectures in
Advanced Mathematics, European Mathematical Society (EMS), Zuerich,
2008.

\bibitem{faS} M. Farber and D. Sch\"{u}tz, Homology of planar polygon
spaces. Geom. Dedicata, 125 (2007),  75-92.



\bibitem{Hatcher} A. Hatcher,  Algebraic topology. Cambridge University
Press, 2002.
%J.-Cl. Hausmann, A. Knutson, \emph{The cohomology ring of polygon
%spaces}. Ann. Inst. Fourier (Grenoble) 48 (1998), 281-321.

%\bibitem{HausmannKnuston2}J.-Cl. Hausmann, A. Knutson,
%\emph{POLYGON SPACES AND GRASSMANNIANS }
\bibitem{Gelfand} I. Gelfand,  M. Goresky, R. MacPherson, and V.
Serganova, Combinatorial geometries, grassmannians, and the moment
map.
 Advances in Math, 63 (1987), 301-316.

\bibitem{Gor} I. Gorodetskaya, Moduli Spaces of Planar Pentagonal Linkages:
Combinatorial Description. arXiv:1305.6756


\bibitem{Kapovich} M. Kapovich,  personal communications.

\bibitem{KapovichMillson} M. Kapovich and J. Millson, On the moduli space
of polygons in the Euclidean plane. J. Diff. Geom., 42 (1995),
430-464.

\bibitem{klya} A. Klyachko, Spatial polygons and stable
configurations of points in the projective line. Tikhomirov,
Alexander (ed.) et al., Algebraic geometry and its applications,
Proceedings of the 8th algebraic geometry conference, Yaroslavl',
Russia, August 10-14, 1992. Braunschweig: Vieweg. Aspects Math. E
25, 67-84 (1994).

\bibitem{PuKho}A.  Pukhlikov and A.  Khovanskii, Finitely additive measures of virtual polyhedra. St. Petersburg Math. J., 4, 2 (1993),
337-356.


\bibitem{PanCycloperm} G. Panina, Cyclopermutohedron,  Trudy Mian,   288
(2015), 149-162.

\bibitem{ziegler}G. Ziegler, Lectures on polytopes. Graduate Texts in Mathematics, vol. 152, Springer-Verlag, New York, 1995.

\bibitem{zvon}D. Zvonkine, Configuration spaces of hinge constructions. Russian
J. of Math. Phys., 5,  2(1997),   247-266.


\end{thebibliography}
\end{document}